# Accelerate Solving Expensive Scheduling by Leveraging Economical Auxiliary Tasks

Minshuo Li, Bo Liu, Bin Xin, Liang Feng, and Peng Li

*Abstract*—To fully leverage the multi-task optimization paradigm for accelerating the solution of expensive scheduling problems, this study has effectively tackled three vital concerns. The primary issue is identifying auxiliary tasks that closely resemble the original expensive task. We suggested a sampling strategy based on job importance, creating a compact matrix by extracting crucial rows from the entire problem specification matrix of the expensive task. This matrix serves as an economical auxiliary task. Mathematically, we proved that this economical auxiliary task bears similarity to its corresponding expensive task and preserves the essential behavior of the expensive task. The subsequent concern revolves around making auxiliary tasks more cost-effective. We determined the sampling proportions for the entire problem specification matrix through factorial design experiments, resulting in a more compact auxiliary task. With a reduced search space and shorter function evaluation time, it can rapidly furnish high-quality transferable information for the primary task. The last aspect involves designing transferable deep information from auxiliary tasks. We regarded the job priorities in the (sub-) optimal solutions to the economical auxiliary task as transferable invariants. By adopting a partial solution patching strategy, we augmented specificity knowledge onto the common knowledge to adapt to the target expensive task. The strategies devised for constructing task pairs and facilitating knowledge transfer, when incorporated into various evolutionary multitasking algorithms, were utilized to address expensive instances of permutation flow shop scheduling. Extensive experiments and statistical comparisons have validated that, with the collaborative synergy of these strategies, the performance of evolutionary multitasking algorithms is significantly enhanced in handling expensive scheduling tasks.

## I. Introduction

EXPENSIVE scheduling is extremely challenging in both fields of combinatorial optimization and evolutionary computation, since they usually have huge search spaces and encounter expensive-to-evaluate objective functions, resulting in time-consuming optimization. Even for the permutation flowshop scheduling problem (PFSP) [1], the type often encountered in practice, when the number of jobs doubles from 100 to 200, its time complexity for evaluating the objective function doubles and the search space becomes $8.45 \times 10^{216}$ times larger, resulting in an unacceptable increase in the cost of exhaustive searches. Consequently, new algorithms that efficiently solve expensive scheduling are urgently needed.

Conventional single-task optimization algorithms, solving problems or tasks in isolation, are stuck in bottlenecks when solving expensive scheduling problems. Combinatorial explosion in the solution space prevalent in expensive scheduling renders the enumeration-based exact methods infeasible. Meanwhile, expensive-to-evaluate objective function in expensive scheduling makes approximation algorithms that rely on iterative search extremely time-consuming, hindering the discovery of excellent solutions from a huge space in a short time. These single-task algorithms, whether exact or approximate, solve different tasks separately and do not take advantage of commonalities between tasks, no matter how similar they are. If commonalities between tasks can be found and leveraged appropriately, it may speed up solving expensive scheduling tasks.

Multi-task optimization, which addresses multiple tasks simultaneously by leveraging their similarities, has paved a promising way for efficiently tackling expensive scheduling, albeit with challenges. It is well known that the effectiveness of multitasking algorithms is highly sensitive to the similarity between tasks [2]. The higher the similarity, the stronger the common knowledge, making multitasking algorithms that utilize common knowledge transfer more effective; and vice versa. Therefore, to take full advantages of the multitasking optimization paradigm and accelerate solving expensive scheduling problems, three crucial concerns need to be addressed. The first concern (*Q1*) is how to identify auxiliary tasks that closely resemble the original expensive task. The second concern (*Q2*) is how to make the auxiliary tasks more cost-effective, as economical auxiliary tasks are easier to solve. The last concern (*Q3*) is how to transform the commonalities between the economical auxiliary tasks and the expensive original task into transferable knowledge, thereby achieving an accelerated search.

Among the aforementioned three concerns, the core lies in how to construct auxiliary tasks that are relatively close to the primary one. In precedent studies, there are two approaches to

This work was supported by the National Key R&D Program of China under grant 2023YFA1009300. *(Corresponding author: Bo Liu)*

M. Li and B. Liu are with Academy of Mathematics and Systems Science, Chinese Academy of Sciences, Beijing 100190 (e-mail: liminshuo@amss.ac.cn (ML), bliu@amss.ac.cn (BL)).

B. Xin is with School of Automation, Beijing Institute of Technology, Beijing 100081, China (e-mail: brucebin@bit.edu.cn).

L. Feng is with College of Computer Science, Chongqing University, Chongqing 400044, China (e-mail: liangf@cqu.edu.cn).

P. Li was with Academy of Mathematics and Systems Science, Chinese Academy of Sciences, Beijing. He is now with the Artificial Intelligence Department, Cainiao Network Group, Hangzhou (e-mail: lipeng@amss.ac.cn).

This article has supplementary material at http://ieeexplore.ieee.org.

constructing task pairs, namely, selecting auxiliary tasks for the primary one. One approach is random selection. Among the ground-breaking multitasking algorithms, evolutionary multitasking (EMT) algorithms randomly selected auxiliary tasks among test benchmarks [2, 3], while multi-task Bayesian optimization randomly extracted subsets from the entire dataset as auxiliary tasks [4]. Another approach involves generating auxiliary tasks based on rules, mainly focusing on continuous optimization [5, 6], while research in combinatorial optimization is rare. In a handful of studies, *K*-means divided all the vertices in a vehicle routing instance into several groups, with the vertices within each group serving as an auxiliary task [7]. Overall, in previous studies, it remains unclear to what extent the constructed auxiliary tasks are related to the primary task. Indeed, measuring the distance between tasks is highly challenging [8, 9]. However, to our knowledge, auxiliary tasks have not been applied to expensive scheduling problems. Identifying closely related auxiliary tasks for a given expensive scheduling problem—the primary task—remains challenging.

Furthermore, if commonalities among tasks are identified, leveraging these commonalities appropriately may expedite problem-solving. Regarding the transfer of deep information from auxiliary tasks, research in the multi-task setting involves both implicit and explicit knowledge transfer. Implicit knowledge transfer was achieved through genetic chromosome crossover [2, 3], while explicit knowledge transfer was achieved through denoising-autoencoder-based solution mapping [10] or by perturbing the best solutions found so far [11]. In fact, how best to learn transferable deeper knowledge—such as partial solutions or dead-end knowledge—remains an open question [12]. To foster search, multi-modality knowledge was designed, where explicit knowledge (partial solutions and complete solutions) and novel implicit knowledge (solution evolution) were exploited and exchanged [9]. In combinatorial optimization, the work in [9] is the first step towards transferring multi-modality knowledge. It does not yet exit, to our knowledge, transferring common knowledge obtained from economical auxiliary tasks, such as partial solutions, to original expensive scheduling task.

In response to the three aforementioned concerns (*Q1-3*), this study aims to construct economical auxiliary tasks closely related to the original expensive scheduling task, design transferable deep commonalities from auxiliary tasks, and fully leverage the multi-task optimization paradigm to expedite the resolution of the primary expensive task. Our contributions are as follows.

1) We introduced a series of sampling strategies based on job importance. By extracting rows corresponding to important jobs from the whole problem specification matrix of the expensive task, a compact matrix is formed, serving as an economical auxiliary task. Among these sampling strategies, the job-importance measure based on largest sum of squares of processing time is most effective in accurately identifying critical jobs in the expensive task, preserving only the essential rows that significantly influence the expensive task's behavior. We have mathematically proven that the economical auxiliary task obtained through this strategy is closely associated with its corresponding expensive task. Moreover, we utilized factorial design to determine the sampling ratio for the entire problem specification matrix, striving to ensure that the economical auxiliary task was as concise as possible. This more concise and closely related auxiliary task to the original task has a smaller search space and shorter evaluation time for the objective function. This enables it to rapidly provide high-quality transferable information for the primary task.

2) To identify and leverage transferable commonalities among tasks, on one hand, we characterized the optimal or suboptimal solutions of economical auxiliary tasks as common knowledge among tasks, and prioritize among jobs as transferable invariances. On the other hand, we introduced a recursive insertion-based strategy for patching partial solutions, appending specificity knowledge onto common knowledge to adapt to the target expensive task. Through the knowledge transfer based on explicit partial solution patching, explicit common knowledge (partial solutions) obtained from economical auxiliary tasks is transferred to the original expensive scheduling task to accelerate its convergence.

3) We integrated the economical auxiliary task (Contribution 1) and the knowledge transfer based on explicit partial solution patching (Contribution 2) into several well-established evolutionary multitasking algorithms to assess their efficacy. We selected expensive instances from the permutation flowshop scheduling problem (PFSP) [1], an extensively studied problem in the literature known for its notorious intractability. Comprehensive numerical experiments and statistical comparisons confirmed that, with the collaborative synergy of both strategies, evolutionary multitask algorithms yielded superior solutions and accelerated convergence when tackling expensive scheduling tasks.

## II. PRELIMINARY

This section introduces the mathematical formulation of the permutation flowshop scheduling problem (PFSP) [13]. Expensive instances derived from this problem are employed to evaluation the effectiveness of the proposed strategies. It also introduces an inter-task distance metric to measure the similarity between permutation flow shop scheduling instances [9]. Based on this inter-task distance metric, we subsequently demonstrate the similarity between economical auxiliary tasks and their expensive counterpart.

### A. Permutation Flowshop Scheduling Problem (PFSP)

PFSP finds a permutation, say a sequence of jobs to be processed on machines, with respect to certain objective(s). In PFSP, a set of $n$ jobs $\{1, 2, \cdots, n\}$ has to be processed on each of the $m$ machines $\{1, 2, \cdots, m\}$. Each machine can execute at most one job at a time, and each job can be executed on at most one machine. The permutation is kept the same on each machine. The $n$ jobs' permutation is denoted as $\pi = [\pi_1, \pi_2, \cdots, \pi_n]$, where the $i$-th element $\pi_i, i \in \{1, 2, \cdots, n\}$ is the job in the $i$-th position of the permutation. The processing time of job $\pi_i$ on machine $j$ is given as $p_{\pi_i,j}$. The completion time for job $\pi_i$ on machine $j$ is denoted as $C_{\pi_i,j}$. The objective is to find a permutation $\pi$ to minimize the maximum completion time for all jobs on all machines, i.e., $C_{\pi_n,m}$. The maximum completion time (makespan) is computed via the recursive (2)-(4) [13].

$$C_{\pi_1,1} = p_{\pi_1,1} \quad (1)$$
$$C_{\pi_i,1} = C_{\pi_{i-1},1} + p_{\pi_i,1}, i = 2, \ldots, n \quad (2)$$
$$C_{\pi_1,j} = C_{\pi_1,j-1} + p_{\pi_1,j}, j = 2, \ldots, m \quad (3)$$
$$C_{(\pi_i,j)} = \max[C_{\pi_{i-1},j}, C_{\pi_i,j-1}] + p_{\pi_i,j}, i,j \geq 2 \quad (4)$$
$$f(\pi) = C_{\pi_n,m}. \quad (5)$$

*B. Inter-Task Distance Metric between PFSPs*

This section introduces a normalized, symmetrical inter-task distance metric, which quantitatively measures the similarity between different PFSP instances. We briefly provide definitions, theorems, and calculations here. Interested readers may refer to [9].

*1) Scale- and Shift- Invariance for PFSP with Makespan*:

**Definition 1.** Problems $f$ and $f'$ are order-isomorphic when $f(\pi) \leq f(\pi') \Leftrightarrow f'(\pi) \leq f'(\pi')$ for any solutions $\pi$ and.

**Theorem 1.** For any two problem specification matrices $P$ and $P'$, and for any positive scale value $t > 0$, if $P' = t \cdot P$, then for an arbitrary solution $\pi$, $f_{P'}(\pi) = t \cdot f_P(\pi)$ is true, where the problem specification matrix $P$ is the $(n \times m)$ matrix with $p_{i,j}$ as its element, $p_{i,j}$ represents the processing time of the $i$-th job on the $j$-th machine. $f_P(\pi)$ and $f_{P'}(\pi)$ denote the makespans with $P$ and $P'$ under solution $\pi$, respectively.

**Theorem 2.** For any two problem specification matrices $P$ and $P'$, If $P' = P + b \cdot E$, where $b \cdot E$ is a shift matrix, $b \in R$, $R$ is the real number set and $E$ is an $(n \times m)$ matrix with all elements equal to 1, then $f_{P'}(\pi) = f_P(\pi) + (m + n - 1) \cdot b$ holds for any arbitrary solution $\pi$.

Theorems 1 and 2 stated that PFSP possesses scale- and shift- invariance property with respect to makespan, respectively. By Definition 1, $f_{P'}$ and $f_P$ are order-isomorphic.

Based on the Theorems 1 and 2, the set for order-isomorphic problems can be defined.

**Definition 2.** An order-isomorphic problem set for function $f_P$ is defined as $G_p = \{f_{P'} | P' = t \cdot P + b \cdot E, \ t > 0, \ b \in R\}$.

By Definition 2, the order-isomorphic problem set $G_p$ contains all problems generated by performing scaled or/and shifted operations on $P$. Geometrically, $P$ can be represented as a point in $R^{n \times m}$ space, while $G_p$ is a set of rays in $R^{n \times m}$ space that are parallel to the ray $\overrightarrow{OP}$ and have different intercepts, where $O$ is the origin of the space.

*2) Inter-task Distance Metric*: The distance $d(f_Q, f_P)$ from $f_Q$ to $f_P$ was defined as the minimum distance from $f_Q$ to the order-isomorphic problems set $G_P$ of $f_P$, and regard the distance as residual error (i.e., the difference) between matrix $Q$ and its optimal approximation $P'$ in the space $G_P$, which is formulated as constrained quadratic programming of

$$d(f_Q, f_P) = \min_{f_{P'} \in G_P} \|Q - P'\|_F = \min_{t>0,b} \|Q - t \cdot P - b \cdot E\|_F. \quad (6)$$

where $\|\cdot\|_F$ is the Frobenius norm making Problem (6) convex, $f_{P'}$ is one function from the order-isomorphic problems set $G_P$, $P'$ is the problem specification matrix of $f_{P'}$, and $t$, $b$ and $E$ were defined in Theorems 1 and 2.

The optimum $t^*$ can be explicitly represented as
$$t^* = \max[t^0, 0], \quad (7)$$
$$t^0 = \frac{n \cdot m \cdot \sum_{i=1}^{n} \sum_{j=1}^{m} q_{i,j} \cdot p_{i,j} - \left(\sum_{i=1}^{n} \sum_{j=1}^{m} q_{i,j}\right) \cdot \left(\sum_{i=1}^{n} \sum_{j=1}^{m} p_{i,j}\right)}{n \cdot m \cdot \sum_{i=1}^{n} \sum_{j=1}^{m} q_{i,j}^2 - \left(\sum_{i=1}^{n} \sum_{j=1}^{m} q_{i,j}\right)^2}. \quad (8)$$

The optimum $b^*$ is
$$b^* = \frac{1}{n \cdot m} \left[\sum_{i=1}^{n} \sum_{j=1}^{m} q_{i,j} - t^* \cdot \left(\sum_{i=1}^{n} \sum_{j=1}^{m} p_{i,j}\right)\right], \quad (9)$$
where $q_{i,j}$ and $p_{i,j}$ are elements of $Q$ and $P$, respectively.

By substituting $t^*$ and $b^*$ into (6), the distance is
$$d(f_Q, f_P) = \|Q - t^* \cdot P - b^* \cdot E\|_F. \quad (10)$$

Then, an analytical distance can be obtained as:
$$d(f_Q, f_P) = \|Q^* - t^* \cdot P^*\|_F, \quad (11)$$
where
$$Q^* = Q - \frac{1}{n \cdot m} \cdot \sum_{i=1}^{n} \sum_{j=1}^{m} q_{i,j} \cdot E. \quad (12)$$
$$P^* = P - \frac{1}{n \cdot m} \cdot \sum_{i=1}^{n} \sum_{j=1}^{m} p_{i,j} \cdot E. \quad (13)$$

*3) Normalization of the inter-task distance*: In accordance with the Cauchy-Buniakowsky-Schwarz inequality, the following inequality holds
$$\|Q^* - t^* \cdot P^*\|_F^2 \geq (\|Q^*\|_F - t^*\|P^*\|_F)^2. \quad (14)$$

Then a normalized distance is
$$d(f_Q, f_P) = \begin{cases} \frac{\|Q^*\|_F - t^*\|P^*\|_F}{\|Q^* - t^* \cdot P^*\|_F}, & \|Q^* - t^* P^*\|_F \neq 0 \\ 0, & otherwise \end{cases}. \quad (15)$$

The normalization of inter-task distance can be illustrated geometrically in Fig. 1. $Q^*$ in (12) and $P^*$ in (13) can be depicted as two points in $R^{n \times m}$ space. Geometrically, (11) can be interpreted as the minimum distance from point $Q^*$ to the ray $tP^*$ $(t \geq 0)$. $\theta$ is the angle between the ray $tP^*(t \geq 0)$ and the ray $tQ^*(t \geq 0)$.

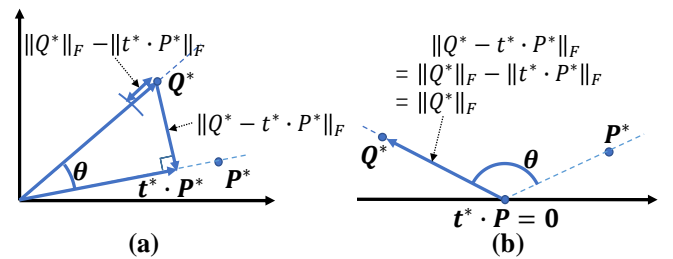

**Fig. 1.** Geometric interpretation of the inter-task distance when $\theta$ is less than $\pi/2$ (**a**) or greater than or equal to $\pi/2$ (**b**) [9].

When $\theta < \pi/2$, as shown in Fig.1(a), we have:
$$\sin \theta = \|Q^* - t^* \cdot P^*\|_F / \|Q^*\|_F, \quad (16)$$
$$\cos \theta = \|t^* \cdot P^*\|_F / \|Q^*\|_F. \quad (17)$$

Then the normalized distance (15) can be represented as:
$$d(f_Q, f_P) = (1 - \cos \theta) / \sin \theta. \quad (18)$$

Obviously, within the interval $[0, \pi/2]$, as $\theta$ increases, the distance function (18) gradually increases, with its value rising from 0 to 1. A distance value of 1 indicates no similarity between PFSPs, while a value of 0 implies that the PFSPs are identical in terms of the order-isomorphic relationship.

When $\theta \geq \pi/2$, as shown in Fig.1(b), the minimum point in the ray $tP^*(t \geq 0)$ to the point $Q^*$ is the origin, i.e., $t^* = 0$. In this regard, the inter-task distance $\|Q^* - t^* \cdot P^*\|_F$ is equal to $\|Q^*\|_F$, and $d(f_Q, f_P) = 1$ according to (15).

Next, we will verify the distance between the economical auxiliary task and the original expensive task based on the geometrically represented distance metric.

### III. CONSTRUCTION OF ECONOMICAL AUXILIARY TASKS AND FOSTERING OF KNOWLEDGE TRANSFER

To fully leverage the advantages of the multitasking optimization paradigm and accelerate solving expensive scheduling problems, this section focuses on addressing two key concerns raised in the research. The first is how to identify auxiliary tasks that are most similar to the original expensive task (Q1). The second is how to transform the commonalities between economical auxiliary tasks and the expensive original task into transferable knowledge for accelerating the search process (Q3).

#### A. Job-Importance based Sampling Strategy to Construct Economical Auxiliary Tasks

It is generally believed that the behavior of scheduling problems is determined by a subset of critical jobs [14-16]. If these critical jobs can be accurately identified, then the economical auxiliary tasks containing them can exhibit behavior that closely aligns with the original expensive task. However, accurately identifying the critical jobs in expensive scheduling tasks remains challenging.

Next, we will provide a measure of job importance and mathematically prove that the economical auxiliary task containing only these most important jobs is closely related to their expensive counterpart.

*1) Job-Importance Measure based on Largest Sum of Squares of Processing Time (LSP)*: We design a method to measure the importance among jobs and sorted the jobs according to their importance values. Specifically, we define the sum of the squares of the processing time of a job on all machines as the importance measure of the job. The higher the value on the metric, the higher the importance of the job; and vice versa. The policy is named as the largest sum of squares of processing time, abbreviated as LSP. We then use a percentage $k$ to select jobs that rank in the top $k$ percent in terms of importance, thereby controlling the size of the subset of critical jobs. We only extract the rows corresponding to the subset of critical jobs from the problem specification matrix of the expensive scheduling task, in order to form a more compact matrix, which constitutes the economical auxiliary task (EAT). In the end, in Section V, we will use factorial design to determine the sampling ratio (the value for $k$) of the entire problem specification matrix, to ensure that EAT is as parsimonious as possible.

Next, we give an example of constructing EAT using LSP-based importance measure, as shown in Fig. 2. The leftmost part of Fig. 2 shows the problem specification matrix of an expensive scheduling task (a), including 10 jobs (from J1 to J10) and 5 machines (from M1 to M5). This task is not expensive and is for illustrative purposes only. Elements in this matrix represent processing time. For example, the element "71" in the fourth row and first column represents the processing time of job J4 on machine M1. Next, we calculate the sum of the squares of each job's processing time on all machines to obtain the job's importance value (b). Each job is ranked according to its importance value. The higher the importance value, the higher the ranking; vice versa. Among these ten jobs, job J4 has the highest importance value and is ranked first (c). We set the value of $k$ to 40, indicating that the top 40% of jobs ranked in importance are selected. In this example, jobs J4, J5, J7 and J9 are selected. The rows corresponding to these four jobs are extracted from the problem specification matrix of the expensive task to form the economical auxiliary task (d).

| M1 | M2 | M3 | M4 | M5 | | | | | | | | | | | | | | |
|----|----|----|----|----|----|----|----|----|----|----|----|----|----|----|----|----|----|----|
| 54 | 79 | 16 | 66 | 58 | J1 | 17133 | J1 | J4 | | | | | | | | | | |
| 83 | 3 | 89 | 58 | 56 | J2 | 21319 | J2 | J9 | | | | | | | | | | |
| 15 | 11 | 49 | 31 | 20 | J3 | 4108 | J3 | J5 | | M1 | M2 | M3 | M4 | M5 | | | | |
| 71 | 99 | 15 | 68 | 85 | J4 | 26916 | J4 | J7 | | 71 | 99 | 15 | 68 | 85 | J4 | | | |
| 77 | 56 | 89 | 78 | 53 | J5 → | 25879 | J5 → | J2 | → | 77 | 56 | 89 | 78 | 53 | J5 | | | |
| 36 | 70 | 45 | 91 | 35 | J6 | 17727 | J6 | J8 | | 87 | 56 | 64 | 85 | 13 | J7 | | | |
| 87 | 56 | 64 | 85 | 13 | J7 | 22195 | J7 | J10 | | 87 | 86 | 75 | 77 | 18 | J9 | | | |
| 76 | 3 | 7 | 85 | 86 | J8 | 20455 | J8 | J6 | | | | | | | | | | |
| 87 | 86 | 75 | 77 | 18 | J9 | 26843 | J9 | J1 | | | | | | | | | | |
| 68 | 5 | 77 | 51 | 68 | J10 | 17803 | J10 | J3 | | | | | | | | | | |

(a) Original task    (b) Measure jobs' importance    (c) Rank jobs based on importance    (d) Construct EAT

**Fig. 2.** Illustrative example of building EAT using LSP-based importance measure.

*2) Mathematical Proof of the Distance between EAT and Its Expensive Counterpart Task*: This section will use the preparatory knowledge of Section II-B to give a mathematical proof that EAT is closely related to its expensive counterpart.

Given that $P = (p_{i,j})_{n \times m}$ is the problem specification matrix for the expensive task. Let $S = \{s_1, s_2, \dots, s_g\}$ be the subset of critical jobs selected from the $n$ jobs of the expensive task, where $s_i \in \{1, \dots, n\}$, $i \in \{1, \dots, g\}$, $g = \lfloor n \cdot k\% \rfloor$, and $\lfloor \cdot \rfloor$ denotes the floor function. Here, $k$ is the sampling ratio defined in Section III-A. Once $S$ is determined, EAT is specified. It is known from Section III-A that $g$ is less than $n$, indicating that the EAT is smaller in scale than the original expensive task. To make the mathematical proof for the distance measure between tasks of the same size in Section II-B applicable, we make the size of EAT consistent with that of the original task by appending rows with elements of 0 to the problem specification matrix of EAT. The zero-padding operation for the problem matrix $Q = (q_{i,j})_{n \times m}$ of EAT is specifically shown in (19).

$$q_{i,j} = \begin{cases} p_{i,j}, & i \in S, \ j = 1, \cdots, m \\ 0, & i \notin S, \ j = 1, \cdots, m \end{cases}. \tag{19}$$

It is easy to prove that EAT is completely equivalent before and after zero padding, that is, makespan is not affected.

In this way, the distance between EAT and its corresponding expensive task can be defined using (6) in Section II-B. Let $A(P)$ and $A(Q)$ denote the sum of all elements in matrices $P$ and $Q$, respectively,

$$A(P) = \sum_{i=1}^{n} \sum_{j=1}^{m} p_{i,j}, \tag{20}$$

$$A(Q) = \sum_{i=1}^{n} \sum_{j=1}^{m} q_{i,j} = \sum_{\alpha \in S} \sum_{j=1}^{m} p_{\alpha,j}. \tag{21}$$

We adopt the symbols defined in Section II-B, where $p_{i,j}^*$ and $q_{i,j}^*$ represent the elements of $P^*$ and $Q^*$, and according to (12) and (13), we have:

$$p_{i,j}^* = p_{i,j} - \frac{1}{nm} A(P), i = 1, \cdots, n, j = 1, \cdots, m, \tag{22}$$

$$q_{i,j}^* = \begin{cases} p_{i,j} - \frac{1}{nm} A(Q), & i \in S, j = 1, \cdots, m \\ -\frac{1}{nm} A(Q), & i \notin S, j = 1, \cdots, m \end{cases}. \tag{23}$$

The normalized distance (18) can be reformulated as:

$$d(f_Q, f_P) = \frac{1-\cos\theta}{\sin\theta} = \sqrt{\frac{2}{1+\cos\theta} - 1}. \tag{24}$$

Minimizing (24) is equivalent to maximizing $\cos\theta$.

According to the geometric interpretation of the inter-task distance, $\theta$ is the angle between the rays $tP^*(t \geq 0)$ and $tQ^*(t \geq 0)$, so

$$\cos\theta = (P^* \cdot Q^*)/(\|P^*\|_2 \cdot \|Q^*\|_2)$$

$$= \frac{\sum_{i=1}^{n} \sum_{j=1}^{m} p_{i,j}^* q_{i,j}^*}{\sqrt{\sum_{i=1}^{n} \sum_{j=1}^{m} (p_{i,j}^*)^2} \cdot \sqrt{\sum_{i=1}^{n} \sum_{j=1}^{m} (q_{i,j}^*)^2}}. \tag{25}$$

We estimate the lower bound of (25) using the following lemma and theorems. Subsequently, we demonstrate that EAT obtained by the job-importance measure based on LSP can maximize this lower bound.

**Lemma 1.** Let $I_1$ and $I_2$ be two sets of positive integers, where $I_2 \subset I_1$. Assuming that $I_1$ and $I_2$ have $N_1$ and $N_2$ elements respectively, then for any $b_{i_1} \in R$, $i_1 \in I_1$, the following inequality holds

$$\left(\sum_{i_1 \in I_1} b_{i_1}\right) \cdot \left(\sum_{i_2 \in I_2} b_{i_2}\right) \leq \frac{N_2}{2} \sum_{i_1 \in I_1} b_{i_1}^2 + \frac{N_1}{2} \sum_{i_2 \in I_2} b_{i_2}^2. \tag{26}$$

*Proof.*

$$\left(\sum_{i_1 \in I_1} b_{i_1}\right) \cdot \left(\sum_{i_2 \in I_2} b_{i_2}\right) = \sum_{i_1 \in I_1} \sum_{i_2 \in I_2} b_{i_1} \cdot b_{i_2}. \tag{27}$$

According to the arithmetic-geometric average inequality [17], for any $i_1 \in I_1$ and $i_2 \in I_2$, the following inequality holds

$$b_{i_1} \cdot b_{i_2} \leq (b_{i_1}^2 + b_{i_2}^2)/2. \tag{28}$$

Thus

$$\left(\sum_{i_1 \in I_1} b_{i_1}\right) \cdot \left(\sum_{i_2 \in I_2} b_{i_2}\right) \leq \sum_{i_1 \in I_1} \sum_{i_2 \in I_2} (b_{i_1}^2 + b_{i_2}^2)/2$$
$$= \frac{N_2}{2} \sum_{i_1 \in I_1} b_{i_1}^2 + \frac{N_1}{2} \sum_{i_2 \in I_2} b_{i_2}^2 \tag{29}$$

∎

**Theorem 3.** $\sum_{i=1}^{n} \sum_{j=1}^{m} p_{i,j}^* q_{i,j}^* \geq \frac{1}{2}(\|Q\|_F^2 - \frac{g}{n}\|P\|_F^2)$ holds.

*Proof.*
$$\sum_{i=1}^{n} \sum_{j=1}^{m} p_{i,j}^* q_{i,j}^* = \sum_{i \in S} \sum_{j=1}^{m} p_{i,j}^* q_{i,j}^* + \sum_{i \notin S} \sum_{j=1}^{m} p_{i,j}^* q_{i,j}^*. \tag{30}$$

By (22) and (23), when $i \in S$, we have

$$p_{i,j}^* q_{i,j}^*$$
$$= \left(p_{i,j} - \frac{1}{nm}A(P)\right)\left(p_{i,j} - \frac{1}{nm}A(Q)\right) \tag{31}$$
$$= p_{i,j}^2 - \frac{1}{nm}p_{i,j}(A(P) + A(Q)) + \frac{1}{n^2m^2}A(P)A(Q),$$

and when $i \notin S$, we have

$$p_{i,j}^* q_{i,j}^* = \left(p_{i,j} - \frac{1}{nm}A(P)\right)\left(-\frac{1}{nm}A(Q)\right)$$
$$= -\frac{1}{nm}p_{i,j}A(Q) + \frac{1}{n^2m^2}A(P)A(Q). \tag{32}$$

By combining (31) and (32), we have

$$\sum_{i=1}^{n}\sum_{j=1}^{m} p_{i,j}^* q_{i,j}^* = \sum_{i \in S}\sum_{j=1}^{m} p_{i,j}^2 - \frac{1}{nm}A(P)A(Q)$$
$$= \|Q\|_F^2 - \frac{1}{nm}A(P)A(Q). \tag{33}$$

According to Lemma 1, (20) and (21),

$$A(P)A(Q) \leq \frac{gm}{2}\sum_{i=1}^{n}\sum_{j=1}^{m} p_{i,j}^2 + \frac{nm}{2}\sum_{i \in S}\sum_{j=1}^{m} p_{i,j}^2$$
$$= \frac{gm}{2}\|P\|_F^2 + \frac{nm}{2}\|Q\|_F^2. \tag{34}$$

Thus

$$\sum_{i=1}^{n}\sum_{j=1}^{m} p_{i,j}^* q_{i,j}^* \geq \|Q\|_F^2 - \frac{1}{mn}\left(\frac{gm}{2}\|P\|_F^2 + \frac{nm}{2}\|Q\|_F^2\right)$$
$$= \frac{1}{2}\left(\|Q\|_F^2 - \frac{g}{n}\|P\|_F^2\right). \tag{35}$$

∎

**Theorem 4.** $\sqrt{\sum_{i=1}^{n}\sum_{j=1}^{m}(p_{i,j}^*)^2}\sqrt{\sum_{i=1}^{n}\sum_{j=1}^{m}(q_{i,j}^*)^2} \leq \left(1 - \frac{1}{mn}\right)\|P\|_F^2$ holds.

*Proof.* We have

$$\sum_{i=1}^{n}\sum_{j=1}^{m}(p_{i,j}^*)^2$$
$$= \sum_{i=1}^{n}\sum_{j=1}^{m}(p_{i,j} - \frac{1}{mn}A(P))^2$$
$$= \sum_{i=1}^{n}\sum_{j=1}^{m}(p_{i,j}^2 - \frac{2}{mn}p_{i,j}A(P) + \frac{1}{m^2n^2}A(P)^2) \tag{36}$$
$$= \|P\|_F^2 - \frac{1}{mn}A(P)^2$$
$$\leq \|P\|_F^2 - \frac{1}{mn}\|P\|_F^2$$
$$= \left(1 - \frac{1}{mn}\right)\|P\|_F^2.$$

Similarly,

$$\sum_{i=1}^{n}\sum_{j=1}^{m}(q_{i,j}^*)^2 \leq \left(1 - \frac{1}{mn}\right)\|Q\|_F^2. \tag{37}$$

Thus,

$$\sqrt{\sum_{i=1}^{n}\sum_{j=1}^{m}(p_{i,j}^*)^2}\sqrt{\sum_{i=1}^{n}\sum_{j=1}^{m}(q_{i,j}^*)^2}$$
$$\leq \sqrt{\left(1-\frac{1}{mn}\right)\|P\|_F^2}\sqrt{\left(1-\frac{1}{mn}\right)\|Q\|_F^2} \tag{38}$$
$$= \left(1-\frac{1}{mn}\right)\|P\|_F\|Q\|_F$$
$$\leq \left(1-\frac{1}{mn}\right)\|P\|_F^2.$$

∎

**Theorem 5.** The lower bound of $\cos\theta$ is determined by $\|Q\|_F^2$.

*Proof.* By substituting (35) and (38) into (25), we obtain

$$\cos\theta \geq \frac{1}{2}\left(\|Q\|_F^2 - \frac{g}{n}\|P\|_F^2\right)/\left(1-\frac{1}{mn}\right)\|P\|_F^2$$
$$= \frac{m}{2(nm-1)}\left(n\frac{\|Q\|_F^2}{\|P\|_F^2} - g\right). \tag{39}$$

As $n$, $m$, $g$ and $P$ are pre-determined, the lower bound of $\cos\theta$ is actually determined by $\|Q\|_F^2$. The larger $\|Q\|_F^2$ is, the larger the lower bound is. ∎

Let $LSP_i$ be the importance value of job $i$, we have

$$\|Q\|_F^2 = \sum_{i \in S}\sum_{j=1}^{m} p_{i,j}^2 = \sum_{i \in S} LSP_i, S = \{s_1, s_2 \ldots, s_g\}. \tag{40}$$

Since the number of selected jobs $g$ is given in advance, it can be known from (40) that in order to maximize $\|Q\|_F^2$, the first $g$ jobs with the highest LSP importance value should be

selected from the complete set of jobs. At this point, we have proven that the EAT constructed by the proposed LSP-importance based sampling strategy can ensure that $\cos\theta$ has an excellent lower bound, thus ensuring the closeness between EAT and the original expensive task.

In Section V, for the purpose of comparison, we designed multiple sampling strategies for generating EAT. We will experimentally compare the similarity between EATs generated by different sampling strategies and their expensive counterparts to elucidate the superiority of the proposed LSP-importance based sampling strategy.

*B. Patching of Partial Solutions to Foster Transferrable Knowledge from the Economical Auxiliary Task to the Original Expensive Task*

After addressing concerns related to generating auxiliary economical tasks that are closest to the original expensive task, another challenge we face is how to leverage the similarity of tasks to accelerate convergence on the original expensive problem. To identify and leverage transferable commonalities among tasks, on the one hand, it is necessary to characterize the common knowledge between tasks, namely, transferable invariant. On the other hand, it is essential to specify the specificity knowledge attached to common knowledge to adapt to the target expensive problem.

*1) Common Knowledge across Different Tasks*: We characterize the optimal or near-optimal solutions of economical auxiliary tasks as common knowledge among tasks, and prioritize among jobs as transferable invariances.

The optimal or near-optimal solutions to EAT, refers to the arrangement of critical jobs selected from the expensive task to meet the objective of minimizing the makespan. Inspired by explicit partial solution methods [9], we choose the optimal solution to EAT as the common knowledge.

The mathematical proof in Section III-A supports the belief behind the aforementioned choice. The behavior of scheduling problems is influenced by certain critical jobs. The job importance measure based on the largest sum of squares of processing time (LSP) effectively identifies critical jobs in expensive scheduling tasks. The EAT, which includes these critical jobs, has been proven to be closest to the original expensive task, exhibiting behavior that is roughly consistent with the original expensive task.

Furthermore, in permutation problems like PFSP, it is highly likely that the precedence between two jobs remains the same in solutions to similar tasks. Therefore, we take the precedence between jobs in the (sub-) optimal solution of the EAT as a transferable invariance, serving as the skeleton for the solution to the original expensive task.

*2) Specificity Knowledge Attached to Common Knowledge*: We propose a recursive insertion-based strategy for patching partial solutions, appending specificity knowledge onto common knowledge to adapt to the target expensive task. Taking (sub-) optimal solution of EAT as the skeleton for the solution to the original expensive task, we recursively insert jobs, which is not selected by job-importance based sampling strategy one by one, eventually the solution of EAT is effectively transformed into a solution for the original expensive task, embodying the application of transferable knowledge in a concrete manner.

In Algorithm 1, we provide a detailed description of the procedures for patching partial solutions based on recursive insertion. $\pi_{EAT}$ represents an optimal solution for EAT, specifically the optimal sequence of critical jobs selected by the job-importance-based sampling strategy to minimize the makespan. $U$ represents the set of alternative jobs not selected as critical jobs. EXP and $\pi_{EXP}$ represent, respectively, the expensive original task and its solution.

Before the recursive insertion begins, $\pi_{EAT}$ serves as the skeleton for $\pi_{EXP}$. When the iteration begins, Algorithm 1 selects the job with the highest importance value from set $U$ and inserts it into all possible positions in the $\pi_{EXP}$ sequence, including the beginning, between any two adjacent jobs, or at the end. This generates a series of candidate partial solutions. In this process, the relative precedence between any pairs of critical jobs forming the skeleton remains unchanged. Evaluate all the above candidate partial solutions and choose the one with the minimum makespan as $\pi_{EXP}$. Remove the selected job from the set $U$. Repeat the above steps until $U$ becomes an empty set, and output the final complete $\pi_{EXP}$.

**Algorithm 1: Patching partial solution based on recursive insertion (RI).**

**Input:** $\pi_{EAT}$: optimal solution to EAT; $U$: the set of jobs not selected for the set of critical jobs; **EXP**: the original expensive task.
**Output**: $\pi_{EXP}$: the complete solution to be transferred to EXP.
1: Set $\pi_{EXP} \coloneqq \pi_{EAT}$.
2: **while** ($U$ is not empty) **do**
3.  Select the job with the largest importance value from $U$.
    *Remark*: The job with larger importance value is preferred to other jobs in $U$, since it brings larger perturbation which is beneficial at the initial stages of constructing solutions.
4.  Insert the selected job into $\pi_{EXP}$ at any possible position, i.e., at the beginning, between any two adjacent jobs, or at the end, to generate a series of candidate solutions.
5.  Evaluate the makespan of all candidate solutions obtained in the previous step.
6.  Select the candidate solution with the smallest makespan as $\pi_{EXP}$.
7.  Remove the selected job from $U$.
8: **end while**

We provide an illustrative example of patching partial solution based on recursive insertion in Fig. 3 to explain the construction process of the solution for EXP. We used EXP and EAT from Fig. 2 as examples. (a) $\pi_{EAT} = [5,9,4,7]$ represents the optimal solution for EAT, and $U = \{2,8,10,6,1,3\}$ represents a set of jobs not selected as critical jobs, arranged in descending order of importance values. Among them, J2 is the job with the highest importance value in $U$ and will be selected first. (b) J2 is inserted into all five possible positions of $\pi_{EAT}$, generating a total of five candidate solutions, and their makespans are calculated. The sequence [5,9,4,2,7] is selected because it has the lowest makespan, and then J2 is removed from $U$. Repeat steps 2-8 of Algorithm 1 until $U$ is an empty set.

This section proposes two strategies aimed at enhancing the capabilities of evolutionary multitasking algorithms in solving expensive scheduling problems. One of the strategies is a sampling strategy based on a job importance measure, used to construct economical auxiliary tasks. The other strategy is a solution patching strategy that converts the optimal solution to the EAT into a solution to the original expensive problem.

To our knowledge, there is currently no research on constructing economical auxiliary tasks for expensive scheduling tasks. Additionally, there is a lack of research on patching solutions to economical auxiliary tasks to address the target expensive task. To compare with the proposed LSP-based job importance measure, we designed seven additional job importance metrics (Section V-A-1). We will demonstrate through experiments that the common knowledge obtained from the EAT based on LSP exhibits superior performance. To compare with the proposed recursive insertion-based partial solution patching strategy, we designed three additional partial solution patching strategies (Section V-A-3). We will demonstrate that the RI-based strategy can more effectively append specificity knowledge on the common knowledge, adapting to the target expensive task, and accelerating its convergence.

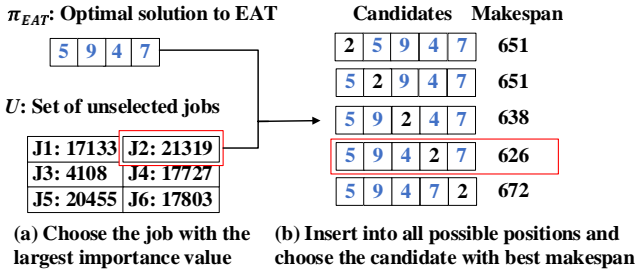

**Fig. 3.** Illustrative example of patching partial solution based on recursive insertion.

## IV. MULTITASKING ALGORITHMS LEVERAGING ECONOMICAL AUXILIARY TASKS TO TACKLE EXPENSIVE SCHEDULING

This section delves into the integration of the two strategies outlined in Section III—economical auxiliary task construction and partial solution patching strategy—into the Evolutionary Multitasking (EMT) algorithms. The aim is to develop multitasking algorithms that utilize economical auxiliary tasks to accelerate the resolution of expensive tasks.

### A. A Brief Introduction to the Evolutionary Multitasking

Evolutionary Multitasking (EMT) utilizes the implicit parallelism of population-based evolutionary search and a genetic information transfer mechanism to concurrently address multiple tasks [2]. This study selected four main EMT algorithms as the carriers for our proposed strategies. We specifically introduced the Multifactorial Evolutionary Algorithm (MFEA-I) [2]. The other three EMT algorithms, namely MFEA-II [3], G-MFEA [5], and P-MFEA [18, 19], are adaptations built upon MFEA-I. For the sake of brevity, we will only outline their distinctions from MFEA-I. Interested readers can refer to the above-mentioned references for algorithm details.

*1) MFEA-I [2]*: It operates as follows: 1) Initialize a population of $N$ individuals. Each individual is represented as a $D$-dimensional real-valued vector using random key encoding, where the elements of the vector take values in the range of 0 to 1. Here, $D$ represents the maximum dimensionality of decision variables for all tasks. 2) Initialize the skill factor for each individual. The skill factor is defined as the identifier of the task in which the individual exhibits relatively higher performance compared to other tasks. After the initialization, the iterative process of genetic evolution commences. 3) Assortative mating and skill factor inheritance. If two randomly selected parents possess the same skill factor or meet a specified random mating probability ($rmp$), they undertake Simulated Binary Crossover (SBX) for offspring reproduction; otherwise, the parents undergo Gaussian mutation. Then, the offspring inherits its parents' skill factor if parents have identical skill factor or inherits an arbitrary parent's skill factor if parents have different skill factors. The offspring is only evaluated on the task that matches its skill factor. 4) Individual learning. Improve each individual using local search operators. 5) Population updating. Use steady-state replacement (i.e., the $\mu + \lambda$ principle) and elitism strategy to update the population. Select the best $N$ individuals from a mixed pool of $\mu$ parents (here, $\mu$ equals $N$) and $\lambda$ offsprings. The evaluation of individual fitness is based on $\varphi_i = 1/\min_{j \in \{1,\cdots,T\}} r_j^i$, where $T$ is the total number of tasks and $r_j^i$ is the rank of individual $i$ on task $j$. Repeat steps (3) to (5) until the termination criteria are met. It is noteworthy that the random mating probability ($rmp$), used to regulates the mating behavior between individuals with different skill factors, is a key parameter for controlling the extent of implicit knowledge transfer across tasks. A value close to 0 for $rmp$ implies that crossover only occurs between parents with the same skill factor, while a value close to 1 allows for completely random mating between parents with different skill factors. In MFEA-I, $rmp$ is fixed at 0.3.

*2) MFEA-II [3]*: To minimize the tendency of harmful knowledge transfer introduced by the fixed $rmp$ in MFEA-I, MFEA-II learned an adaptive transfer parameter matrix to guide the transfer intensity across tasks at runtime, replacing the predetermined $rmp$. The transfer parameter matrix is purely data-driven, learned by minimizing the Kullback-Leibler (KL) divergence between the probability distributions of offspring populations and parent populations across all tasks. The transfer parameter matrix serves as a surrogate for the similarity between different tasks, facilitating adaptive knowledge transfer across related tasks.

*3) G-MFEA [5]*: By incorporating two strategies, namely decision variable translation and shuffling strategy, into MFEA-I, the generalized MFEA (G-MFEA) was formulated. To mitigate the impact of differences in optimal solutions for weakly correlated tasks on the performance of MFEA-I, the decision variable translation set the optimal solutions for all tasks to be the same. To address tasks with different dimensions, the shuffling strategy randomly reordered the sequence of decision variables and allowed each variable to transfer knowledge across tasks. By leveraging high-quality solutions from multiple computationally cheap tasks, G-MFEA effectively reduced the number of evaluations on expensive tasks, solving expensive continuous optimization efficiently.

*4) P-MFEA [18, 19]*: To address permutation-based combinatorial optimization problems, P-MFEA replaced the real-valued encoding in MFEA-I with permutation-based encoding. Additionally, it employed ordered crossover and swap mutation to replace SBX and Gaussian mutation in MFEA-I, respectively.

*B. Incorporation of the Economical Auxiliary Task and Partial Solution Patching Strategy into EMT*

A multi-task algorithm has been developed for addressing expensive scheduling tasks, as illustrated in Algorithm 2. This is achieved by integrating the design of economical auxiliary tasks (Step 1 highlighted in bold) and the strategy for patching partial solutions (Steps 10 to 16 highlighted in bold) into the standard EMT algorithms. It is highly anticipated that these enhancements will expedite the solution of expensive tasks by effectively leveraging transferable knowledge from economical auxiliary tasks.

*1) A Pair of Tasks Composed of the Primary Expensive Task (EXP) and its Economical Auxiliary Task (EAT)*: In multi-task settings, a pair of tasks is solved simultaneously. In EMT, task pairs are often randomly selected from benchmark [2, 3]. In our study, we construct a closely related economical auxiliary task (EAT) for a given expensive scheduling task (primary task), as detailed in Section III-A. Thus, the EXP and the EAT constitute a task pair, which is simultaneously solved in a multi-task setting.

*2) Solutions to EAT, after Being Repaired by the Partial Solution Patching Strategy, are Transferred to EXP:* In addition to preserving the default implicit knowledge transfer based on assortative mating in EMT, the best solutions from the EAT in the current generation are patched into complete solutions for the EXP through the recursive insertion-based partial solution patching strategy (detailed in Section III.B and Algorithm 1). This set of patched solutions serves as transferrable explicit knowledge, migrating to form new high-quality solutions for the EXP. We implement the settings for transfer triggering conditions and the number of transferred individuals as described in [20]. Throughout the iterative process of genetic evolution, every $G$ generations (with $G$ set to 5), we select the top-performing $S$ individuals (with $S$ set to 5) possessing EAT skill factors from the current population, perform patching, and then transfer them.

**Algorithm 2: Multitasking algorithm leveraging transferable knowledge from economical auxiliary tasks to tackle expensive task.**

**Input:** EXP, the targeted expensive PFSP task.
**Output:** $\pi^*_{EXP}$, the best solution obtained for the targeted expensive PFSP task.

1: **Construct a pair of tasks consisting of the primary expensive scheduling task (EXP) and its economical auxiliary task (EAT).**
   *Remark*: We employ the job-importance based sampling strategy to create a similar economical auxiliary task for a given expensive scheduling task (primary task), as described in Section III-A and Section V-A-1.
2: Initialize a population of $N$ individuals.
3: Initialize the skill factor $\tau_i$ for each individual $i$, $(i = 1, ..., N)$.
   *Remark*: If individual $i$ performs better on the expensive task (or economical auxiliary task), the value of $\tau_i$ is the string "EXP" (or "EAT").
4: Set $gen = 1$.
5: **While** (the termination criterion is not satisfied) **do**
6:    Apply assortative mating and skill factor inheritance to generate the offspring population.
7:    Decode the offspring individuals from the real-valued vectors using Ranked-Order Value rule into permutations, as detailed in Section IV-C-1.
8:    Improve the permutation corresponding to each individual using local search, as detailed in Section IV-C-2.
9:    Adjust the real-valued vector to correspond to the improved permutation, as detailed in Section IV-C-3.

   // Solutions to the economical auxiliary task, after being repaired by the partial solution patching strategy, are migrated to the original expensive task.

10:    **If ($\mathrm{mod}(gen, G) == 0$)**
11:      **Select the top-performing $S$ individuals among those with the skill factor "EAT" in current population.**
12:      **Use the Ranked-Order Value rule to decode the selected $S$ individuals into permutations for the economical auxiliary task, denoted as $\pi_{EAT}$, as described in Section IV-C-1.**
13:      **Use the partial solution patching strategy to patch $\pi_{EAT}$ into a complete solution $\pi$ for the expensive task, as detailed in Section III-B and Algorithm 1.**
14:      **Retrieve the real-valued encoding $x$ corresponding to the permutation $\pi$, as detailed in Section IV-C-4.**
15:      **The $x$ is added as a new individual to the offspring population and assigned the skill factor "EXP", as detailed in Section IV-C-5.**
16:    **End If**
17:    Concatenate current population and offspring population, and calculate the fitness value for each individual using the following equation:
$$\varphi_i = 1/\min_{j \in \{1, \cdots, T\}} r^i_j$$
where $T$ is the total number of tasks and $r^i_j$ is the rank of individual $i$ on task $j$.
18:    Use steady-state replacement and elitism strategy to update the population.
19:    $gen = gen + 1$.
20: **End While**
21: Select the individual with the highest fitness among those with the "EXP" skill factor in the population, convert it into a job sequence using Ranked-Order Value rule. This job sequence represents the best solution achievable for the targeted expensive PFSP task, denoted as $\pi^*_{EXP}$.

*C. Miscellaneous Items*

After introducing the construction of economical auxiliary tasks and the strategy of patching partial solutions into EMT, we present additional details to consider when solving expensive PFSP.

*1) Decoding into Scheduling Solutions*: The four aforementioned EMT methods employ two different encoding schemes, namely job-permutation encoding and real-valued vector encoding. The encoding based on permutation itself constitutes the scheduling solution. Here, a method is provided for mapping individuals encoded with real numbers into permutations.

MFEA-I, MFEA-II, and G-MFEA utilize a *D*-dimensional real-valued vector encoding based on random key, where $D$ represents the number of jobs in the original expensive scheduling task. We employ the Ranked-Order Value (ROV) rule [13] to decode the real-valued vector $x = [x_1, \cdots, x_D]$ into a feasible schedule, that is, a permutation $\pi = [\pi_1, \cdots, \pi_D]$. Hereafter, for simplicity, we use $\pi$ and $x$ to respectively represent the scheduling solution of the original expensive scheduling task and its corresponding real-valued vector. According to the ROV rule, $\pi_l$ represents the ranking of $x_l$ in the real-valued vector when sorted in ascending order. To obtain the permutation $\pi_{EAT}$ for the economical auxiliary tasks, jobs not belonging to EAT are removed from $\pi$.

Next, we will demonstrate how to map an individual encoded with a real-valued vector into a permutation. Suppose the real-valued vector $x$ is [0.61, 0.65, 0.01, 0.86, 0.97, 0.69, 0.99, 0.63, 0.78, 0.29]. According to the ROV rule, $x$ is decoded into the permutation $\pi$ for EXP, which is [3, 5, 1, 8, 9, 6, 10, 4, 7, 2]. Taking EAT in Fig. 2 as an example, since EAT only includes jobs 4, 5, 7, and 9, removing other jobs from $\pi$ results in the permutation for EAT, $\pi_{EAT}$ = [5, 9, 4, 7].

*2) Individual Improvement*: The local search methods employed by the four EMT approaches mentioned above cannot be directly applied to PFSP. MFEA-I, MFEA-II, and G-MFEA do not address scheduling, while the N6 neighborhood-based local search in P-MFEA is applicable to job shop but not suitable for PFSP. Thus, we utilize the INSERT-based local search [13] to enhance the performance of individuals. For fair comparison, the original local searches in the four standard EMT algorithms have been replaced with the INSERT-based local search.

It operates as follows: randomly select two distinct jobs from the permutation $\pi$, and then insert the latter job before the former one. For $\pi_{EAT}$, select any two jobs belonging to the EAT, and insert the latter job before the former one. For the permutation of each offspring, the INSERT-based local search is performed with a search intensity of $L$ iterations, and the best solution $\pi^{ls}$ obtained during this process is selected.

*3) Adjusting the Solution Encoded with a Real-valued Vector*: After implementing individual improvement to the permutation, it is necessary to adjust the individual encoded as a real-value vector to ensure its alignment with the improved permutation. Since the ROV rule is employed when converting real-valued vectors into permutations, achieving mutual correspondence between the two is straightforward. After identifying jobs whose positions have changed in the improved permutation, reposition the real-number elements corresponding to those jobs in the real-valued vector to ensure their alignment with the positions of those jobs in the new permutation. Given that P-MFEA directly operates on permutations, the above operation is unnecessary. Below, we provide an example to explain the process of adjusting the solution encoded with a real-valued vector.

We take the real-valued vector $x$ and its corresponding permutation $\pi$ from Miscellaneous Item (1) as an example. Suppose an individual improvement is applied to the permutation $\pi$ = [3, 5, 1, 8, 9, 6, 10, 4, 7, 2], resulting in the improved permutation $\pi^{ls}$ =[1, 3, 5, 8, 9, 6, 10, 4, 7, 2]. It can be observed that the positions of jobs 1, 3, and 5 have changed. We adjust the positions of the corresponding real-valued elements in vector $x$, resulting in the improved vector $x^{ls}$ = [0.01, 0.61, 0.65, 0.86, 0.97, 0.69, 0.99, 0.63, 0.78, 0.29].

*4) Map the Transferred Permutation to a Real-valued Vector*: In the context of MFEA-I, MFEA-II, and G-MFEA, after patching $\pi_{EAT}$ with a partial solution patching strategy into a complete solution $\pi$ for the expensive task, it is necessary to obtain the real-valued encoded solution $x$ corresponding to $\pi$. We obtain the real-valued vector corresponding to this permutation using the method described in Miscellaneous Item (3).

*5) Assigning Skill Factors to Transferred Solutions*: After mapping the transferred permutation to a real-valued vector, the individual encoded with real numbers is added to the offspring population and assigned the skill factor "EXP".

*6) Learning of the Adaptive Transfer Parameter Matrix in MFEA-II*: It is not affected by the introduction of the economical auxiliary tasks and the partial solution patching strategy into MFEA-II.

*7) About the Two Strategies in G-MFEA*: The two strategies in G-MFEA—decision variable translation and shuffling strategy—are unaffected and retained.

*8) Time for Constructing EAT*: It is included in the algorithm's runtime.

## V. EXPERIMENTS TO FIND THE BEST STRATEGIES FOR CONSTRUCTING ECONOMICAL AUXILIARY TASKS AND TRANSFERRING KNOWLEDGE

This section is designed to analyze experiments with the goal of determining the best strategy combination for application in the context of EMT. This includes: i) constructing economical auxiliary tasks (Q1), ii) determining the sampling ratio influencing the size of economical auxiliary tasks (Q2), and iii) patching partial solutions (Q3).

### A. Factors and Levels for Comparison

We design the levels corresponding to the three factors as follows.

*1) Job-importance based Sampling Strategies to Construct Economical Auxiliary Tasks*: To compare with the proposed job importance measure based on LSP, we designed six additional job importance measures derived from the well-known Nawaz-Enscore-Ham (NEH) heuristic [21] and its variants [22, 23]. Additionally, we also designed a random sampling strategy.

LSP (Largest Sum of Squares of Processing Time): The importance value of a job is defined as the sum of the squares of its processing time across all machines. See Section III-A for details.

LST (Largest Sum of Processing Time): The importance value of a job is defined as the sum of its processing time across all machines.

KK1: Inspired by NEHKK1 [22], the importance value of job $i$ is defined as $\min\{a_i, b_i\}$, where $a_i = \sum_{j=1}^{m}[(m-1)(m-2)/2 + m - i] \cdot p_{i,j}$ and $b_i = \sum_{j=1}^{m}[(m-1)(m-2)/2 + i - 1] \cdot p_{i,j}$.

KK2: Inspired by NEHKK2 [23], the importance value of job $i$ is defined as $\min\{a_i, b_i\}$, where $a_i = T_i + U_i$, $b_i = T_i - U_i$, $T_i = \sum_{j=1}^{m} p_{i,j}$, $U_i = \sum_{j=1}^{m}\frac{j - 3/4}{\lfloor m/2 \rfloor - 3/4} \cdot (p_{i,\lfloor m/2 \rfloor + 1 - j} - p_{i,\lceil m/2 \rceil + j})$, $\lfloor \cdot \rfloor$ and $\lceil \cdot \rceil$ denote the floor and ceil functions, respectively.

SR0: The importance value of a job is defined as its position in the permutation obtained by NEH.

SR1: The importance value of a job is defined as its position in the permutation obtained by NEHKK1.

SR2: The importance value of a job is defined as its position in the permutation obtained by NEHKK2.

RND: The importance value of a job is defined as its position in a randomly generated permutation.

Under the first four measures (LSP, LST, KK1, and KK2), a higher value indicates higher job importance; and vice versa. In the latter four measures (SR0, SR1, SR2, and RND), the closer a job is positioned to the beginning of the sequence, the higher its importance; and vice versa.

*2) Sampling Ratio of the Entire Problem Specification Matrix (Percentage $k$)*: We introduced a percentage parameter $k$ to select jobs based on their importance values, choosing those ranked in the top $k$ percent. This parameter determines the size of the critical job set, specifies the number of rows extracted from the problem specification matrix of the expensive scheduling task, thereby controlling the size of the economical auxiliary tasks. The values of $k$ are set to be evenly spaced from 10 to 90, with an interval of 10.

*3) Partial Solution Patching Strategies*: To compare with the proposed partial solution patching strategy based on recursive insertion (RI), we additionally designed three alternative partial solution patching strategies. Similar to RI, these three strategies select the job with the highest importance value from $U$ (representing the set of jobs not selected as critical jobs, arranged in descending order of importance), insert it into the (near-) optimal solution to EAT, and then remove the job from $U$. The process iterates until the partial solution is patched into a complete solution for the expensive task. The difference among these strategies lies in the different insertion positions of the jobs.

RI (Recursive Insertion): See Algorithm 1 for details.

EI (Insert at End): Inspired by [24], insert the job at the end of the current partial permutation.

OI (Odd/even dependent Insertion): Inspired by [25], if the current permutation's length is odd, insert the job at the end of the permutation; otherwise, insert the job at the beginning of the permutation.

AI (Arbitrary Insertion): Insert the job at an arbitrary position in the permutation.

*B. Benchmark*

We selected expensive instances from Taillard's benchmark [26]. The benchmark consists of 12 groups of instances with different sizes, where the number of jobs ($n$) takes values $\{20, 50, 100, 200, 500\}$, and the number of machines ($m$) takes values $\{5, 10, 20\}$. For each combination of job quantity and machine quantity, there are 10 instances, resulting in a total of 120 instances. The computational cost for optimizing instances where $n \times m \geq 500$ is quite high [27]. There is a total of 80 instances that meet this condition, namely, ta41-50 (50 jobs, 10 machines), ta51-60 (50 jobs, 20 machines), ta61-70 (100 jobs, 5 machines), ta71-80 (100 jobs, 10 machines), ta81-90 (100 jobs, 20 machines), ta91-100 (200 jobs, 10 machines), ta101-110 (200 jobs, 20 machines), and ta111-120 (500 jobs, 20 machines). In this study, we selected these 80 instances as the original expensive tasks.

*C. Performance Metrics*

*1) Inter-task Distance Metric*: We employed the Inter-task Distance Metric (ITDM) from (15) to assess the similarity between economical auxiliary tasks (EAT) generated by different sampling strategies and the original expensive task (EXP). The smaller the value of ITDM, the more similar the tasks are; conversely, a larger value indicates greater dissimilarity.

*2) Algorithm Performance*: Average Relative Error (ARE), Best Relative Error (BRE), and Worst Relative Error (WRE) are used, refer to (41-43). $RE_{i,j,l}$ represents the relative error of algorithm $i$ on instance $j$ in the $l$-th run (44). $C_{i,j,l}$ is the makespan obtained by algorithm $i$ on instance $j$ in the $l$-th run, and $C_j^*$ is the best-so-far solution for instance $j$. $L$ denotes the number of independent runs, set to 20.

$$ARE_{i,j} = (1/L) \cdot \sum_{l=1}^{L} RE_{i,j,l} \qquad (41)$$
$$BRE_{i,j} = \min_{l=1,\cdots,L}\{RE_{i,j,l}\} \qquad (42)$$
$$WRE_{i,j} = \max_{l=1,\cdots,L}\{RE_{i,j,l}\} \qquad (43)$$
$$RE_{i,j,l} = 100 \times (C_{i,j,l} - C_j^*)/C_j^* \qquad (44)$$

*D. Comparison of Distances between EXP and EATs Generated by Different Importance Sampling Strategies*

This sub-section evaluates the effectiveness of various strategies in constructing economical auxiliary tasks. A total of 72 strategies were formed by combining 8 sampling strategies based on different job importance measures with 9 sampling ratios. On the 80 expensive instances from Taillard's benchmark, we calculated the distance between the EATs constructed by different strategies and the original expensive task (EXP). The box plot in Fig. 4 illustrates, at various sampling ratios, the average distance values between EAT generated by different sampling strategies and the original EXP across all instances. From Fig. 4, it is evident that, at each sampling ratio, the LSP-based strategy produces an EAT that is the closest to the original expensive task when compared to other importance sampling strategies.

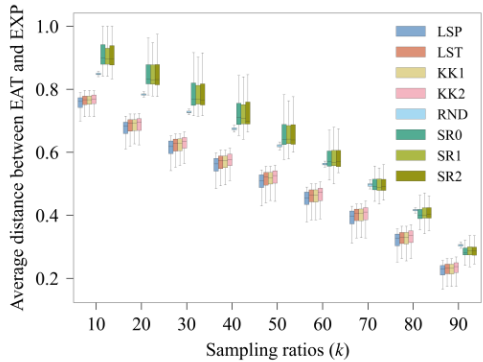

**Fig. 4.** Box plots of the average distance between EAT and EXP across all instances for different sampling strategies and sampling ratios.

We further examined the performance of each importance sampling strategy across all sampling ratios, as shown in Fig. 5. Wilcoxon signed-rank test [28] at a 95% confidence level was employed to assess the statistical differences between the results. At a significance level of 0.05, LSP was significantly superior to all other importance sampling strategies: LST ($p = 5.79 \times 10^{-115}$), KK1 ($p = 1.88 \times 10^{-115}$), KK2 ($p = 4.62 \times 10^{-118}$), RND ($p < 10^{-118}$), SR0 ($p < 10^{-118}$), SR1 ($p < 10^{-118}$) and SR2 ($p < 10^{-118}$). Additionally, LSP, LST, KK1, and KK2 outperformed the remaining four strategies.

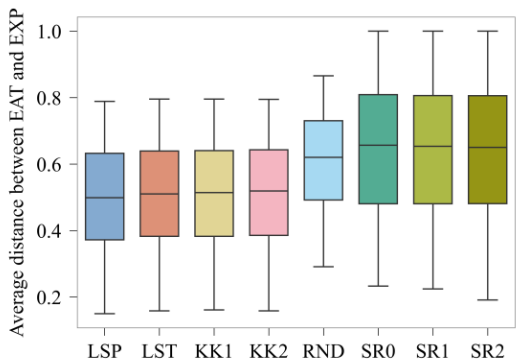

**Fig. 5.** Box plots of the average distance between EAT and EXP across all instances and sampling ratios for different importance sampling strategies.

*E. Comparison of Partial Solution Patching Strategies*

This sub-section examines the effectiveness of the four partial solution patching strategies. For each expensive instance, we can obtain a total of 72 EATs by using 8 importance sampling strategies and 9 sampling ratios. We had a total of 80 expensive tasks, resulting in 5760 EATs. To obtain the (near-) optimal solution for each EAT, we first used the NEH heuristic [21] to find a high-quality initial guess. Subsequently, the initial solution was improved through 10,000 iterations using the simulated annealing as described in [9]. The solution to EAT was patched into a complete solution to EXP using each of the four partial solution patching strategies. ARE was obtained by evaluating their makespan.

Fig. 6 presents box plots of the ARE obtained under different sampling ratios and various partial solution patching strategies across all instances and all importance sampling strategies. It is evident that the partial solution patching strategy based on recursive insertion (RI) outperforms other patching strategies at each sampling ratio.

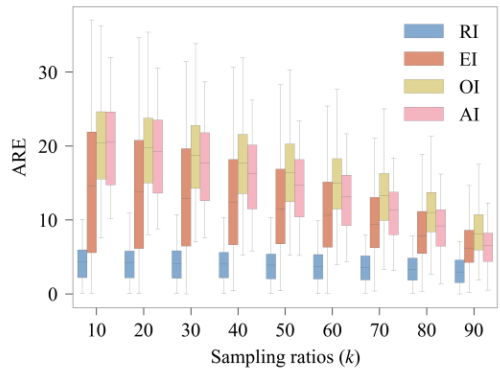

**Fig. 6.** Box plots of the ARE obtained under different sampling ratios and various partial solution patching strategies across all instances and all importance sampling strategies.

Fig. 7 displays box plots of the ARE obtained under different partial solution patching strategies across all instances, importance sampling strategies, and sampling ratios. We conducted the Wilcoxon signed-rank test at a 95% confidence level to examine the statistical differences between the results. Statistical comparisons indicate that RI is significantly superior to EI, OI, and AI, with $p$-values all less than $10^{-118}$.

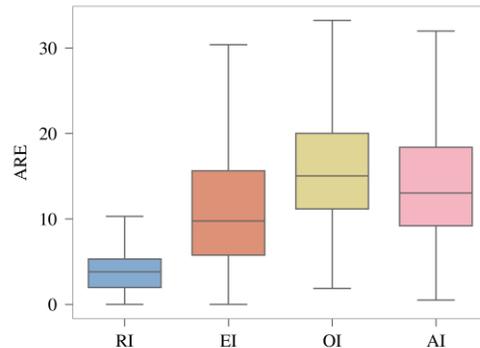

**Fig. 7.** Box plots of the ARE obtained under different partial solution patching strategies across all instances, importance sampling strategies, and sampling ratios.

*F. Identifying Optimal Strategies for Constructing EAT in the EMT Context*

This sub-section conducts a comprehensive factorial design experiment to determine the optimal strategies for constructing EAT in the EMT context. Given the effectiveness of RI compared to the other three partial solution patching strategies, RI was utilized in this experiment. There is a total of 72 strategies for constructing EAT, formed by combining 8 importance sampling strategies with 9 sampling ratios. We incorporated each strategy for constructing EAT and RI-based patching scheme into MFEA-I, the extensively studied EMT algorithm [2], resulting in a total of 72 variants of the MFEA-

I. The expensive instances of the Taillard's benchmark consist of 8 different scales. We randomly selected two instances from each of the 10 instances for each scale, resulting in a total of 16 instances used for testing. The 72 variants of the MFEA-I algorithm were repeated 20 times on each instance. The algorithm termination criterion was set to a maximum CPU runtime, $T = 0.03nm$ seconds, where $n$ and $m$ represent the number of jobs and machines in the instance, respectively.

To differentiate among these 72 strategies, we calculated their ARE performances and conducted multiple comparisons using Tukey's test [28] with a confidence level of 0.05. The results are presented in Table I. When grouping the strategies with Tukey's test, different letter labels represent different groups. There are no significant differences within groups of strategies, while significant differences exist between groups. The results in Table I indicate that the 7 strategies belonging to Group A, namely LSP-20 (representing the importance sampling strategy as LSP with a sampling ratio $k$ of 20), LST-20, KK2-20, LSP-30, KK1-20, LST-30, and KK2-30, have the lowest average ARE values. They are considered the best group under Tukey's test, exhibiting superior performance compared to other strategy groups. Next, we choose these seven strategies to construct EAT.

TABLE I
RESULTS OF A FULL FACTORIAL DESIGN TO FIND THE BEST STRATEGIES FOR CONSTRUCTING EAT IN THE EMT CONTEXT

| Config | ARE | Group by Tukey's Test | Config | ARE | Group by Tukey's Test |
|---|---|---|---|---|---|
| LSP-20 | 2.81 | A | SR0-40 | 4.25 | IJKLMNOPQR |
| LST-20 | 2.81 | A | SR2-30 | 4.29 | JKLMNOPQR |
| KK2-20 | 2.87 | A | SR1-30 | 4.34 | JKLMNOPQRS |
| LSP-30 | 2.87 | A | SR2-40 | 4.38 | JKLMNOPQRS |
| KK1-20 | 2.89 | A | SR1-40 | 4.5 | KLMNOPQRS |
| LST-30 | 2.9 | A | SR0-50 | 4.51 | KLMNOPQRS |
| KK2-30 | 2.91 | A | SR2-50 | 4.61 | LMNOPQRS |
| KK1-30 | 2.94 | AB | SR1-50 | 4.7 | MNOPQRS |
| LST-10 | 2.96 | AB | RND-60 | 4.75 | NOPQRST |
| LST-40 | 3.07 | ABC | KK2-70 | 4.84 | OPQRSTU |
| LSP-40 | 3.08 | ABC | SR0-60 | 4.86 | PQRSTU |
| KK1-10 | 3.11 | ABC | LST-70 | 4.91 | PQRSTU |
| KK2-40 | 3.14 | ABCD | KK1-70 | 4.95 | QRSTU |
| KK1-40 | 3.14 | ABC | SR2-60 | 4.95 | QRSTU |
| KK2-10 | 3.23 | ABCDE | LSP-70 | 4.99 | RSTU |
| LSP-10 | 3.26 | ABCDE | SR1-60 | 5.07 | STU |
| KK2-50 | 3.41 | ABCDEF | SR2-70 | 5.5 | TUV |
| RND-20 | 3.43 | ABCDEFG | RND-70 | 5.51 | UVW |
| LSP-50 | 3.43 | ABCDEFG | SR0-70 | 5.53 | UVW |
| LST-50 | 3.43 | ABCDEFG | SR1-70 | 5.58 | UVWX |
| KK1-50 | 3.48 | ABCDEFGH | KK2-80 | 6.13 | VWXY |
| RND-30 | 3.5 | ABCDEFGHI | LST-80 | 6.26 | WXY |
| SR0-10 | 3.69 | BCDEFGHIJ | KK1-80 | 6.3 | XY |
| RND-10 | 3.74 | CDEFGHIJ | LSP-80 | 6.43 | Y |
| RND-40 | 3.75 | CDEFGHIJK | SR1-80 | 6.49 | Y |
| SR0-20 | 3.9 | DEFGHIJKL | SR0-80 | 6.53 | Y |
| KK2-60 | 3.93 | EFGHIJKL | SR2-80 | 6.54 | Y |
| SR2-10 | 3.97 | EFGHIJKLM | RND-80 | 6.69 | Y |
| LST-60 | 4.02 | FGHIJKLMN | SR1-90 | 8.04 | Z |
| KK1-60 | 4.03 | FGHIJKLMN | SR2-90 | 8.07 | Z |
| LSP-60 | 4.03 | FGHIJKLMN | SR0-90 | 8.12 | Z |
| SR1-10 | 4.04 | FGHIJKLMN | KK2-90 | 8.42 | Z |
| SR0-30 | 4.06 | FGHIJKLMN | KK1-90 | 8.49 | Z |
| SR2-20 | 4.09 | FGHIJKLMNO | LST-90 | 8.51 | Z |
| RND-50 | 4.18 | GHIJKLMNOP | RND-90 | 8.59 | Z |
| SR1-20 | 4.23 | HIJKLMNOPQ | LSP-90 | 8.59 | Z |

*G. Discussion on Preliminary Experimental Results*

We have obtained several preliminary conclusions. Firstly, the preliminary experiments confirm that the EAT generated by the LSP-importance-based sampling strategy is closest to the original expensive task under the ITDM distance metric. These experimental results are consistent with the theoretical conclusions presented in Section III-A. The experimental results from the full factorial design indicate that the combined strategy, utilizing LSP-importance sampling and a sampling ratio of 20, has the highest average performance. And RI is the most effective partial solution patching strategy. Additionally, we observed that the proposed importance sampling strategies, LST, KK1, and KK2, demonstrate performance comparable to LSP. This phenomenon suggests that for other types of combinatorial problems where ITDM cannot be applied to calculate distances between tasks, LST, KK1, and KK2 can be utilized to construct economical auxiliary tasks.

Next, we will conduct large-scale experiments to validate the effectiveness of the aforementioned strategy combinations in solving expensive tasks within the framework of EMT.

VI. NUMERICAL EXPERIMENTS AND COMPARISONS

In this section, we designate two strategies employed in EMT algorithms—specifically, randomly generated task pairs and implicit knowledge transfer—as benchmark strategies for comparison. These two strategies, along with successful strategies from the previous section, including seven strategies for constructing economical auxiliary tasks and one strategy for partial solution patching, are combined to form a strategy pool. Various strategy combinations are generated by selecting strategies from this pool. Each of these strategy combinations is incorporated into all four EMT algorithms, and extensive computational experiments are conducted on large-scale expensive instances. Through statistical analysis, our goal is to validate the impact of combining different auxiliary task construction strategies and knowledge transfer strategies on enhancing the performance of EMT algorithms in solving expensive scheduling problems.

*A. Algorithms for Comparisons*

In this study, any algorithm utilized for comparison is instantiated through the following three elements: the manner in which task pairs are constructed, the method of knowledge transfer between tasks, and the EMT algorithm that implements them.

*1) Construction of Task Pairs*: There are two primary approaches to selecting or constructing auxiliary tasks for the primary task. The first involves the random generation of task pairs, as utilized in standard EMT algorithms, where task pairs are randomly chosen from benchmarks. The second approach is the importance-based sampling strategies proposed in this study, aimed at constructing an economical auxiliary task that is similar to the primary task.

Three methods for randomly generating task pairs (RndTsk1-3): Firstly, auxiliary tasks should have the same number of machines as the primary expensive task to ensure similarity. Next, from Taillard's benchmark, we randomly i) RndTsk1: choose instances with the same number of jobs as the primary task as auxiliary tasks; ii) RndTsk2: choose instances with fewer jobs than the primary task as auxiliary tasks; iii) RndTsk3: choose instances with more jobs than the primary task as auxiliary tasks. When selecting auxiliary tasks using the RndTsk3 method, instances ta61-70, ta91-100, and ta111-120 need to be excluded because there are no instances in Taillard's benchmark with the same number of machines but more jobs than them.

Seven strategies for constructing task pairs based on job-importance sampling (ImpTsk): We selected the top-performing 7 job-importance sampling schemes from Section V-F, namely LSP-20, LST-20, KK2-20, LSP-30, KK1-20, LST-30, and KK2-30, to construct an economical auxiliary task for a given expensive scheduling task.

2) *Knowledge transfer between tasks*: One is implicit knowledge transfer used in EMT (represented as IK), and the other is a recursive insert-based partial solution patching strategy (RI).

3) *EMT algorithms*: The four EMT algorithms introduced in Section IV-A are employed, namely MFEA-I, MFEA-II, G-MFEA and P-MFEA.

4) *Algorithm Representation based on Triplets*: We use the terms "EMT/Task pair/Knowledge transfer" to name the algorithm. Specifically, the set EMT={MFEA-I, MFEA-II, G-MFEA, P-MFEA}; nested set Task pair={RndTsk, ImpTsk}, where RndTsk={RndTsk1, RndTsk2, RndTsk3}, ImpTsk={LSP-20, LST-20, KK2-20, LSP-30, KK1-20, LST-30, KK2-30}; and Knowledge transfer={IK, RI}.

Using the triplet notation, it is easy to create algorithms or sets of algorithms with various configurations. It should be noted that {RndTsk1, RndTsk3} cannot be combined with RI, as RI is only applicable when the number of jobs for auxiliary tasks is less than that of the primary task. Thus, by combining ten methods for constructing task pairs, two knowledge transfer approaches, and four EMT algorithms, and excluding invalid combinations, a total of 72 algorithms are formed.

Based on the above triplet notation, "MFEA-I/LSP-20/RI" signifies that the MFEA-I algorithm employs the LSP-based sampling strategy with a sampling ratio of 20 to construct auxiliary tasks. It also utilizes a recursive insert-based partial solution patching strategy to achieve knowledge transfer between tasks. Similarly, "EMT/{RndTsk, ImpTsk}/IK" represents a collection of 40 algorithms, achieved by combining 4 EMT algorithms, ten methods from the {RndTsk, ImpTsk} set, and one implicit knowledge transfer method.

*B. Computational Environment*

We evaluated the performances of all 72 algorithms on 80 expensive instances. Algorithms terminate when reaching a maximum CPU runtime of $T = 0.03nm$ seconds, where $n$ and $m$ represent the number of jobs and machines for the primary task, respectively. Each algorithm was independently run 20 times on each instance. All algorithms were implemented in Python 3.8.5, with the Cython code from [27] being used for calculating the makespan. Experiments were conducted on an Intel Xeon CPU E5-2650 2.20GHz machine with 128 GB memory running Ubuntu 16.04.7 LTS.

*C. Impacts of Economical Auxiliary Tasks and Knowledge Transfer on Accelerating EMT Algorithms*

We selected the default strategies in EMT algorithms—random generation of task pairs and implicit knowledge transfer—as benchmark policies for comparison. We scrutinized the impact of our proposed strategy for constructing task pairs based on job-importance sampling, as well as the knowledge transfer strategy based on partial solution patching, on the performance of EMT algorithms in addressing expensive scheduling tasks. Table IV's first column enumerates configurations used for comparison, including different task pairs and knowledge transfer strategies. The second column presents their ARE, BRE, and WRE values across all instances. The third column conducts Wilcoxon rank-sum test [28] to detect significant differences in ARE. The fourth column provides the effect sizes of Cohen's $d$ [29] to assess the strength of differences between them. Fig. 8 illustrated the convergence performance of the average ARE over time for algorithm sets configured with different task pair strategies and knowledge transfer strategies across all instances.

TABLE II
PERFORMANCES OF ALGORITHMS CONFIGURED WITH DIFFERENT TASK PAIRS AND TRANSFER STRATEGIES

| Algorithm/Algorithm Set | [ARE BRE WRE] | p-value | Effect Size of Cohen's d |
|---|---|---|---|
| **Effects arising from different task pairs** | | | |
| EMT/RndTsk/{IK, RI} vs. EMT/ImpTsk/{IK, RI} | [9.63 8.21 10.84] vs. [6.67 5.88 7.36] | $< 10^{-69}$ | 0.58 |
| EMT/RndTsk/IK vs. EMT/ImpTsk/IK | [11.70 10.10 12.99] vs. [10.56 9.39 11.52] | $< 10^{-10}$ | 0.26 |
| EMT/RndTsk/RI vs. EMT/ImpTsk/RI | [4.17 3.26 5.19] vs. [2.79 2.36 3.20] | $< 10^{-24}$ | 0.76 |
| **Effects arising from different knowledge transfer** | | | |
| EMT/{RndTsk, ImpTsk}/IK vs. EMT/{RndTsk, ImpTsk}/RI | [10.87 9.59 11.92] vs. [2.96 2.47 3.45] | $< 10^{-200}$ | 2.27 |
| EMT/RndTsk/IK vs. EMT/RndTsk/RI | [11.70 10.10 12.99] vs. [4.17 3.26 5.19] | $< 10^{-105}$ | 1.88 |
| EMT/ImpTsk/IK vs. EMT/ImpTsk/RI | [10.56 9.39 11.52] vs. [2.79 2.36 3.20] | $< 10^{-200}$ | 2.37 |
| **Effects arising from different combinations of task pairs and knowledge transfers** | | | |
| EMT/RndTsk/IK vs. EMT/ImpTsk/RI | [11.70 10.10 12.99] vs. [2.79 2.36 3.20] | $< 10^{-200}$ | 3.17 |

We have the following findings. 1) Compared to the strategy of randomly generating task pairs (RndTsk), the task pair construction strategy based on job-importance sampling (ImpTsk) achieved better search quality in the complete set of knowledge transfer strategies ($p$-value: $4.08 \times 10^{-70}$; effect size of Cohen's $d$: 0.58, representing a medium effect) and under the same knowledge transfer strategies (for IK and RI, the $p$-values are $1.53 \times 10^{-11}$ and $4.59 \times 10^{-25}$, with Cohen's $d$ effect sizes of 0.26 (small effect) and 0.76 (medium effect), respectively). Fig. 8 illustrated that under the same knowledge transfer strategy, ImpTsk exhibited faster convergence compared to RndTsk.

2) Compared to implicit knowledge transfer strategy (IK), explicit knowledge transfer strategy (RI) has achieved better

search quality in the complete set of strategies used for constructing task pairs ($p$-value is $< 10^{-200}$; Cohen's $d$ effect size is 2.27, indicating the presence of a huge effect) and under the same strategies for constructing task pairs (for RndTsk and ImpTsk, the $p$-values are $1.46 \times 10^{-106}$ and $< 10^{-200}$, with Cohen's $d$ effect sizes of 1.88 (very large effect) and 2.37 (huge effect), respectively). Fig. 8 also indicated that under the same strategy for constructing task pairs, RI exhibited faster convergence compared to IK.

3) Compared to the combination of RndTsk and IK, the combination of ImpTsk and RI achieved better search quality ($p$-value is $< 10^{-200}$; Cohen's $d$ effect size is 3.17, indicating the presence of a huge effect), with improvements of 76%, 77%, and 75% in terms of ARE, BRE, and WRE, respectively. Fig. 8 illustrated that the combination of ImpTsk and RI achieved the fastest convergence.

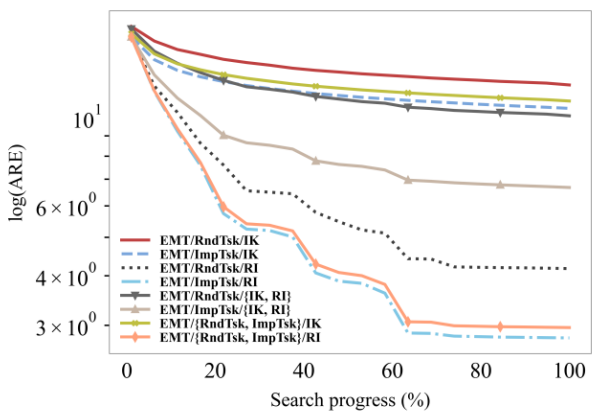

**Fig. 8.** Convergence curves of the average ARE for algorithms configured with different task pairs and knowledge transfer strategies across all instances.

## VII. CONCLUSIONS

Efficiently solving expensive scheduling problems remains challenging, given the expansive search space and the expensive-to-evaluate objective function. This study successfully addressed three key issues: constructing an economical auxiliary task that closely resembles the original expensive task, determining the scale of the economical auxiliary task, as well as identifying the inter-task transferable commonalities, along with the specialized knowledge to adapt to the expensive task. The strategies proposed for constructing task pairs and knowledge transfer have been seamlessly integrated into various evolutionary multitasking algorithms. Comprehensive numerical experiments and statistical comparisons confirmed that, under the combined effect of these strategies, the advantages of the multi-task optimization paradigm are fully triggered and utilized. This resulted in a significant enhancement of the performance of EMT algorithms when addressing expensive scheduling tasks.

We made several assumptions, and limitations exist, both of which should be noted.

In the context of multi-task optimization, the composition of task pairs becomes particularly crucial. The effectiveness of multi-task optimization algorithms is highly dependent on the degree of similarity between tasks. In previous studies, the construction of task pairs—selecting auxiliary tasks for the prime task—was implemented by randomly picking tasks from a set of benchmarks. We propose a job-importance measure based on largest sum of squares of processing time, capable of accurately identifying critical jobs in the expensive task. With the foundation of inter-task distance measure [9], we have mathematically demonstrated that economical auxiliary tasks, containing only these most important jobs, are closely associated with their corresponding expensive task. To the best of our knowledge, this is the first report to identify the closest economical auxiliary tasks for a given expensive scheduling task. It addressed the challenge of quantifying the degree of proximity between the compact auxiliary tasks and the primary task in combinatorial optimization, which has been previously unquantifiable. In this study, we only employed a single economical auxiliary task to enhance the search. A worthwhile avenue for future research is to explore the interaction among multiple economical auxiliary tasks to further enhance search performance. While such an approach may introduce additional computational overhead, it is expected to yield more diverse and transferable knowledge, thereby assisting in further improving search performance on expensive tasks.

*Size of economical auxiliary tasks*: Properly sizing economical auxiliary tasks ensures not only their close relevance to the primary expensive task but also ensures that solving these economical auxiliary tasks comes with a lower cost. We selected jobs ranked in the top $k$ percent in terms of importance as key jobs. Following this, we extracted the rows from the problem specification matrix of the expensive task that correspond to these key jobs, thereby forming the problem specification matrix for auxiliary tasks. Through factorial design experiments, we determined that the number of rows in the problem specification matrix for auxiliary tasks is set to 20% or 30% of the original problem specification matrix's rows. Compared to the expensive primary task, this more compact matrix has a smaller search space and shorter objective function evaluation time, enabling it to rapidly provide transferable information for the primary task. Exploring the balance among the theoretical lower bounds of the scale of economical auxiliary tasks, their similarity to the primary task, and maximizing transferable knowledge is worth in-depth research.

*Identifying and leveraging transferable commonalities among tasks*: The effective transfer of genuine commonalities between tasks will expedite problem-solving. However, how to best learn transferable deeper knowledge—such as partial solutions, deadlock knowledge—remains an unresolved question. On one hand, we characterized the optimal or suboptimal solutions of economical auxiliary tasks as common knowledge among tasks, and prioritize among jobs as transferable invariances. On the other hand, we introduced a recursive insertion-based strategy for patching partial solutions, appending specificity knowledge onto common knowledge to adapt to the target expensive task. To the best of our knowledge, this is the first time that explicit common knowledge (partial solutions) obtained from economical auxiliary tasks has been transferred to the original expensive scheduling task to accelerate its convergence. In this study, there is only unidirectional knowledge transfer from economical auxiliary tasks to the original task. Exploring how to achieve bidirectional knowledge transfer and multimodal knowledge transfer is worth investigating.

In summary, by effectively leveraging the commonalities between the primary expensive task and its closely related economical auxiliary tasks, this study has paved a promising path for efficiently addressing expensive scheduling problems in the context of multi-task optimization.

### REFERENCES

[1] S. M. Johnson, "Optimal two- and three- stage production schedules with setup times included," *Naval Research Logistics Quarterly,* vol. 1, pp. 61-68, 1954.
[2] A. Gupta, Y. S. Ong, and L. Feng, "Multifactorial Evolution: Toward Evolutionary Multitasking," *IEEE Transactions on Evolutionary Computation,* vol. 20, pp. 343-357, 2016.
[3] K. K. Bali, Y. Ong, A. Gupta, and P. S. Tan, "Multifactorial Evolutionary Algorithm With Online Transfer Parameter Estimation: MFEA-II," *IEEE Transactions on Evolutionary Computation,* vol. 24, pp. 69-83, 2020.
[4] K. Swersky, J. Snoek, and R. P. Adams, "Multi-task Bayesian optimization," *Proceedings of the Advances in Neural Information Processing Systems,* pp. 2004–2012, 2013.
[5] J. Ding, C. Yang, Y. Jin, and T. Chai, "Generalized Multitasking for Evolutionary Optimization of Expensive Problems," *IEEE Transactions on Evolutionary Computation,* vol. 23, pp. 44-58, 2019.
[6] Y. Feng, L. Feng, S. Kwong, and K. C. Tan, "A Multivariation Multifactorial Evolutionary Algorithm for Large-Scale Multiobjective Optimization," *IEEE Transactions on Evolutionary Computation,* vol. 26, pp. 248-262, 2022.
[7] Q. Shang, Y. Huang, Y. Wang, M. Li, and L. Feng, "Solving vehicle routing problem by memetic search with evolutionary multitasking," *Memetic Computing,* vol. 14, pp. 31-44, 2022.
[8] Y. Bengio, A. Lodi, and A. Prouvost, "Machine learning for combinatorial optimization: A methodological tour d'horizon," *European Journal of Operational Research,* vol. 290, pp. 405-421, 2021.
[9] P. Li and B. Liu, "Multi-task Combinatorial Optimization: Adaptive Multi-modality Knowledge Transfer by an Explicit Inter-task Distance," Academy of Mathematics and Systems Science, Chinese Academy of Sciences, Beijing, 2020.
[10] L. Feng, L. Zhou, J. Zhong, A. Gupta, Y. S. Ong, K. C. Tan*, et al.*, "Evolutionary Multitasking via Explicit Autoencoding," *IEEE Transactions on Cybernetics,* vol. 49, pp. 3457-3470, 2019.
[11] M. Y. Cheng, A. Gupta, Y. S. Ong, and Z. W. Ni, "Coevolutionary multitasking for concurrent global optimization: With case studies in complex engineering design," *Engineering Applications of Artificial Intelligence,* vol. 64, pp. 13-24, 2017.
[12] S. Toyer, S. Thiebaux, F. Trevizan, and L. X. Xie, "ASNets: Deep Learning for Generalised Planning," *Journal of Artificial Intelligence Research,* vol. 68, pp. 1-68, 2020.
[13] B. Liu, L. Wang, and Y. H. Jin, "An effective PSO-based memetic algorithm for flow shop scheduling," *IEEE Transactions on Systems Man and Cybernetics Part B-Cybernetics,* vol. 37, pp. 18-27, 2007.
[14] J.-P. Watson, L. Barbulescu, L. D. Whitley, and A. E. Howe, "Contrasting Structured and Random Permutation Flow-Shop Scheduling Problems: Search-Space Topology and Algorithm Performance," *INFORMS Journal on Computing,* vol. 14, pp. 98-123, 2002.
[15] J. P. Watson, L. D. Whitley, and A. E. Howe, "Linking search space structure, run-time dynamics, and problem difficulty: A step toward demystifying tabu search," *Journal of Artificial Intelligence Research,* vol. 24, pp. 221-261, 2005.
[16] B. Liu, J.-J. Xu, B. Qian, J.-R. Wang, and Y.-B. Chu, "Probabilistic memetic algorithm for flowshop scheduling," in *IEEE Symposium Series on Computational Intelligence*, 2013, pp. 60-64.
[17] P. S. Bullen, "The Arithmetic, Geometric and Harmonic Means," in *Handbook of Means and Their Inequalities*, P. S. Bullen, Ed., Dordrecht: Springer Netherlands, 2003, pp. 60-174.
[18] Z. Lei, L. Feng, Z. Jinghui, Y. S. Ong, Z. Zhu, and E. Sha, "Evolutionary multitasking in combinatorial search spaces: A case study in capacitated vehicle routing problem," in *2016 IEEE Symposium Series on Computational Intelligence (SSCI)*, 2016, pp. 1-8.
[19] Y. Yuan, Y. S. Ong, A. Gupta, P. S. Tan, and H. Xu, "Evolutionary multitasking in permutation-based combinatorial optimization problems: Realization with TSP, QAP, LOP, and JSP," in *2016 IEEE Region 10 Conference (TENCON)*, 2016, pp. 3157-3164.
[20] L. Feng, Y. Huang, L. Zhou, J. Zhong, A. Gupta, K. Tang*, et al.*, "Explicit Evolutionary Multitasking for Combinatorial Optimization: A Case Study on Capacitated Vehicle Routing Problem," *IEEE Transactions on Cybernetics,* vol. 51, pp. 3143-3156, 2021.
[21] M. Nawaz, E. E. Enscore, and I. Ham, "A heuristic algorithm for the m-machine, n-job flow-shop sequencing problem," *Omega,* vol. 11, pp. 91-95, 1983.
[22] P. J. Kalczynski and J. Kamburowski, "An improved NEH heuristic to minimize makespan in permutation flow shops," *Computers & Operations Research,* vol. 35, pp. 3001-3008, 2008.
[23] P. J. Kalczynski and J. Kamburowski, "An empirical analysis of the optimality rate of flow shop heuristics," *European Journal of Operational Research,* vol. 198, pp. 93-101, 2009.
[24] D. S. Palmer, "Sequencing Jobs through a Multi-Stage Process in the Minimum Total Time - a Quick Method of Obtaining a near Optimum," *Operational Research Quarterly,* vol. 16, pp. 101-107, 1965.
[25] J. N. D. Gupta, "Heuristic Algorithms for Multistage Flowshop Scheduling Problem," *IISE Transactions,* vol. 4, pp. 11-18, 1972.
[26] E. Taillard, "Benchmarks for basic scheduling problems," *European Journal of Operational Research,* vol. 64, pp. 278-285, 1993.
[27] M. Karimi-Mamaghan, M. Mohammadi, B. Pasdeloup, and P. Meyer, "Learning to select operators in meta-heuristics: An integration of Q-learning into the iterated greedy algorithm for the permutation flowshop scheduling problem," *European Journal of Operational Research,* vol. 304, pp. 1296-1330, 2023.
[28] A. M. Mood, F. A. Graybill, and D. C. Boes, *Introduction to the theory of statistics*. New York: McGraw-Hill, 1973.
[29] C. O. Fritz, P. E. Morris, and J. J. Richler, "Effect size estimates: current use, calculations, and interpretation," *Journal of Experimental Psychology: General,* vol. 141, p. 2, 2012.